\documentclass[11pt]{article}
\usepackage{latexsym, amsmath}

\newtheorem{theorem}{Theorem}
\newtheorem{lemma}{Lemma}

\newtheorem{definition}{Definition}

\newtheorem{conjecture}{Conjecture}
\newtheorem{proposition}{Proposition}

\newenvironment{proof}{{\it Proof:\/}}{\hfill $\Box$\\ } 

\newcommand{\Z}{\mbox{\bf Z}}
\newcommand{\F}{\mbox{\bf  F}}

\newcommand{\cc}{\mathbf{c}}

\title{ Primality Proving via One Round in ECPP\\ 
and One Iteration in AKS}

\date{}
 
\author{Qi Cheng\thanks{School of Computer Science,
the University of Oklahoma,
Norman, OK 73019, USA.
Email: {\tt qcheng@cs.ou.edu}
}}

\begin{document}

\maketitle

\begin{abstract}
On August 2002, Agrawal, Kayal and Saxena announced the first
deterministic and polynomial time primality testing algorithm.
For an input $n$, the AKS algorithm runs in heuristical time
$\tilde{O}(\log^6 n)$. Verification takes roughly the same amount
of time.
On the other hand, the Elliptic Curve Primality Proving algorithm (ECPP), 
runs in random heuristical time  $\tilde{O}(\log^{6} n)$
( $\tilde{O}(\log^5 n)$ if the fast multiplication is used), and
generates certificates which can be easily verified.
More recently, Berrizbeitia gave a variant
of the AKS algorithm, in which some primes cost
much less time to prove than a general prime does.
In this paper, we explore the possibility of combining the
ideas in these celebrated algorithms to design a more
efficient algorithm.
A random primality proving algorithm  with 
heuristic time complexity  $\tilde{O}(\log^{4} n)$
is presented. It generates a certificate
of primality which is $O(\log n)$ bits long and 
can be verified in deterministic time $\tilde{O}(\log^4 n)$.
The reduction in  time complexity is achieved by first generalizing 
Berrizbeitia's  algorithm to one which has
higher density of easily-proved primes.
For a general prime, one round of ECPP is deployed to reduce 
its primality proof to the proof of a random easily-proved prime.
\end{abstract}

\section{Introduction}

Testing whether a number is prime or not is one of the fundamental problems
in computational number theory. It has wide applications in computer science, 
especially in cryptography. After tremendous efforts invested by researchers 
in about two hundred years,
it was finally proved by Agrawal, Kayal and Saxena
\cite{AgrawalKa02} that the set of primes is in the complexity class
{\bf P}.  For a given integer $n$, the AKS  algorithm runs in
time no longer than $\tilde{O}( \log^{12} n)$, while the best
deterministic algorithm before it has 
subexponential complexity \cite{AdlemanPo83}.
Under a reasonable conjecture, The AKS algorithm 
should give out answer in time 
$\tilde{O}(\log^6 n)$. 

{\bf Notation:} In this paper, we use ``$\ln$'' for logarithm base $e$ and
``$\log$'' for logarithm base $2$. We write $r^\alpha || n$, if
$r^\alpha | n$ but $r^{\alpha + 1} \not| n$. By $\tilde{O} (f(n))$,
we mean $O( f(n) \mathrm{polylog}(f(n)))$.

The AKS algorithm is based on the derandomization
of a polynomial identity testing. It involves many iterations
of polynomial modular exponentiation.
To test the primality of a integer $n$,
the algorithm first searches for a suitable $r$, which is 
provably $O(\log^6 n)$, or heuristically $O(\log^2 n)$.
Then the algorithm will check for $s$ from $1$ to 
$S= \lceil 2\sqrt{r} \log n \rceil$, 
whether
\begin{equation}\label{aksequation}
(x+s)^n = x^n +s \pmod{n, x^r-1}. 
\end{equation}
The algorithm declares that $n$ is a prime if all the checks
pass. %Under a reasonable conjecture,
%$S =O( (\log n)^2 )$. 
The computing of $(x+s)^n \pmod{n, x^r-1}$ takes time
$\tilde{O}(r \log^2 n)$ if we use the fast multiplication.
The total time complexity is thus $\tilde{O}(rS\log^2 n)$.

While  the AKS algorithm is a great accomplishment in the theory,  
the current version  is very slow. 
Unless its time complexity can be dramatically improved,
it cannot replace random primality testing algorithms
with better efficiency.
In most of applications in cryptography, an efficient random algorithm is 
sufficient, as long as the algorithm can generate
a certificate of primality, which in deterministic time convinces 
a verifier who does not believe any number theory conjectures.
A primality testing algorithm which generates a certificate of
primality is sometimes called {\em primality proving algorithm}.
Similarly a primality testing algorithm which generates a certificate of
compositeness is sometimes called {\em compositeness proving algorithm}.
Very efficient random compositeness 
proving algorithms have long been known.
Curiously, primality proving algorithms 
lag far behind of compositeness proving algorithms in term of
efficiency and simplicity.

Recently, Berrizbeitia \cite{Berrizbeitia02} proposed a brilliant 
modification to the AKS original algorithm. 
He used the polynomial $x^{2^s} - a$ instead of $x^r - 1$
in equation~(\ref{aksequation}), where $2^s \approx \log^2 n$.
Among others, he was able to prove the following proposition:

\begin{proposition}\label{theoremofBerrizbeitia}
Given an integer $n\equiv 1 \pmod{4} $. Denote 
$s = \lceil 2\log\log n\rceil$. Assume that $ 2^k || n-1$ and $k\geq s $.
If there exists an integer $a$, such that $({a \over n}) = -1$
and $a^{ n-1 \over 2} \equiv -1 \pmod{n}$, then
$$ (1 + x)^n \equiv 1 + x^n \pmod{n, x^{2^s}-a} $$
iff $n$ is a power of a prime.
\end{proposition}

Unlike the AKS algorithm, where each prime costs roughly the same,
there are ``easily-proved primes'' in Berrizbeitia's algorithm,
namely, the primes $p$ where $p-1$ has a factor of a power of two
larger than $\log^2 n$.
For those primes, one iteration of polynomial
modular exponentiation, which runs in time $\tilde{O}(\log^4 n)$, 
establishes the primality of $p$, provided that a suitable $a$ exists.
In fact, $a$  can be found easily 
if $n$ is indeed a prime and randomness is allowed
in the algorithm. It serves as a prime certificate for $n$.
\begin{definition}
In this paper, for a primality proving algorithm,
we call a prime $p$ {\em easily-proved}, if the algorithm runs
in expected time $\tilde{O}(\log^4 p)$ on $p$.
\end{definition}

What is the density of the easily-proved primes in Berrizbeitia's 
algorithm? Heuristically for a random prime $p$, $p-1$
should have probability $ 1 \over \log^2 p$ 
to have a factor $2^s \approx \log^2 p$, 
hence the easily-proved primes have density
$1 \over \log^2 p$ around $p$ in his algorithm.

\subsection{Increasing the density of easily-proved primes}

We prove the following theorem in Section~\ref{generalization},
which can be regarded as a generalization of 
Proposition~\ref{theoremofBerrizbeitia}.

\begin{theorem}\label{main}
(Main) Given a number $n$ which is not a power of an integer. 
Suppose that there exists a prime $r$,
$r^\alpha||n-1 ( \alpha \geq  1)$ and 
$ r \geq \log^2 n $.
In addition, there exists a number $1< a < n$, such that
$a^{r^\alpha} \equiv 1 \pmod{n}$,
$gcd(a^{r^{\alpha-1}}-1,n) = 1$, and
$$ ( 1 + x)^n = 1 + x^n \pmod{n, x^r - a}, $$
then $n$ is a prime.
\end{theorem}

The number $a$ can be found easily if $n$ is a prime and
randomness is allowed. It serves as a prime certificate for $n$.
Base on this theorem, we propose
a random algorithm which establishes the primality
of $p$ in time $\tilde{O}(\log^4 p)$ if $p-1$ contains a prime factor 
between $\log^2 p$ and $C \log^2 p$ for some small constant $C$. 
\begin{definition}
We call a positive integer $n$ {\em $C$-good}, if
$n-1$ has a prime factor $p$ such that $\log^2 n \leq p \leq C \log^2 n$.
\end{definition}
What is the density of $C$-good primes?
Apparently the density should be higher than the 
density of easily-proved primes in Berrizbeitia's
algorithm. 
Let $m = \prod_{p\ \mathrm{prime},b_1\leq p\leq C b_1} p.$
First we count the number of integers between $1$ and $m$ which
have a prime factor between $b_1$ and $C b_1$. This is precisely
the number of zero-divisors in ring $\Z/m\Z$:
$$ (m-1) - m \prod_{p\ \mathrm{prime}, 
b_1\leq p\leq C b_1} ( 1 - {1 \over p}). $$
We will prove in Section~\ref{densitysection}
that this number is greater than $ m \over \ln b_1$
for $C=\cc$ and $b_1$ sufficiently large,
where $\cc$ is an absolute constant to be determined later. 
To analyze the time complexity of our algorithm,
we mainly concern the density of $2$-good primes in short intervals.
For simplicity, we call a number {\em good}, when it is $2$-good.
Since compared with $\log^2 n$, $n$ is very big, we expect that 
\begin{conjecture}\label{densityconjecture}
There exists an absolute constant  $\lambda$,
such that for any sufficiently large integer $n$,
$$ {  Number\ of\ 2-good\ primes\ between\ n-2\sqrt{n} +1\ and\ 
n + 2\sqrt{n} +1  
\over Number\ of\ primes\ between\ n-2\sqrt{n} +1\ and\ n + 2\sqrt{n} +1
 } 
> {\lambda \over \ln (\log^2 n)}. $$
\end{conjecture}
We are unable to prove this inequality however, but
we present in the paper some numerical evidences.
We comment that questions about the prime distribution 
in a short interval are usually  very hard to answer.

\subsection{Algorithm for the general primes}

For general primes, we apply the idea in the Elliptic Curve
Primality Proving algorithm (ECPP).
ECPP was proposed by Goldwasser, Kilian \cite{GoldwasserKi86} 
and Atkin \cite{Atkin86} and implemented by 
Atkin and Morain \cite{AtkinMo93}.  In practice, ECPP performs
much better than the current version of AKS. It has been used to
prove primality of numbers up to thousands of decimal 
digits \cite{Morain98}.

In ECPP, if we want to prove that an integer $n$ is a prime,
we reduce the problem to the proof of primality of
a smaller number (less than $n/2$).
To achieve this, we try to find an elliptic curve
with $\omega n'$ points over $\Z/n\Z$, where
$\omega$ is completely factored and $n'$ is a probable prime
greater than $ (\sqrt[4]{n}+1 )^2$. 
Once we have such a curve and a point on the curve with
order $n'$, the primality of $n'$ implies
the primality of $n$.
Since point counting on elliptic curves is  expensive, 
we usually use the elliptic curves
with complex multiplications of small discriminants.
Nonetheless, it is plausible to assume that the order
of the curve has the desired form with the same
probability as a random integer does.
ECPP needs $O(\log n)$ rounds
of reductions to eventually reduce the problem to a primality proof of 
a very small prime, say, less than $1000$. As observed in
\cite{LenstraLe90}, one round of reduction 
takes heuristic time $\tilde{O}(\log^{5} n)$,
or $\tilde{O}(\log^{4} n) $ if we use the fast multiplication.
To get the time complexity, it is assumed that
the number of primes between $n-2\sqrt{n}+1$ and
$n+2\sqrt{n}+1$ is greater than $ \sqrt{n} / \log^2 n $,
and the number of points on an elliptic curve with small
discriminant complex multiplication behaves like 
a random number in the Hassa range. We refer the
assumption as the ECPP heuristics.
Rigorous proof of the time complexity seems out of reach,
as it involves the study of the prime distribution in a short interval.

Our algorithm can be decomposed into two stages.
In the first stage, for a general probable prime $n$, 
we will use one round of ECPP to 
reduce its proof of primality  to a good probable
prime $n'$ near $n$. For convenience,  we require that
$n - 2\sqrt{n} + 1 \leq n' \leq n + 2\sqrt{n} + 1$ (See 
section~\ref{implementation} for implementation issues).
Note that up to a constant factor,
the time complexity of one round reduction in ECPP
is equivalent to the time complexity of finding a curve
with a prime order.
In the set of primes 
between $n -2 \sqrt{n} + 1$ and $n +2 \sqrt{n} + 1$,
the density of good primes
is $\lambda \over \ln (\log^2 n)$ by conjecture.
Hence heuristically {\em the extra condition on $n'$ (that $n'$ should
be good) will increase the time complexity merely by a factor 
of $O(\log \log n)$}. Therefore for all the primes, 
without significant increase of time complexity,
we reduce its primality proving to the proof of a good prime.
In the second stage, we find a primality certificate
for $n'$. To do this, we search for $a$ which satisfies the
conditions in the main theorem, and compute the
polynomial modular exponentiation.
Heuristically, the total expected running time of 
the first and the second stages becomes $\tilde{O}(\log^{4} n) $.
However, due to the short interval of the number of points
over elliptic curves, it seems difficult to obtain the
rigorous time complexity. 
Put it altogether, we now have a general purpose prime proving algorithm,
which has following properties:
\begin{enumerate}
\item it runs very fast ($\tilde{O}( \log^4 n )$ ) assuming 
reasonable heuristics.
\item For many primes, ECPP subroutine is not needed.
\item The certificate, which consists of the curve, a point
on the curve with order $n'$, $n'$ and $a$, is very short. 
It consists of only $O(\log n)$ bits as opposed to $O(\log^2 n)$ 
bits in ECPP.
\item A verifier can be convinced in deterministic
time $\tilde{O}(\log^4 n)$. In fact, the most time consuming 
part in the verification is the iteration of polynomial modular 
exponentiation. 
\end{enumerate}

This paper is organized as following:
In Section~\ref{AKSECPP}, we review the propositions
used by AKS and ECPP to prove primality. 
In Section~\ref{descriptionofalgorithm},
we describe our algorithm and present the time
complexity analysis. In Section~\ref{densitysection},
we prove a theorem which can be regarded as an 
evidence for the density heuristics.
The main theorem is proved in Section~\ref{generalization}.
We conclude this paper with some discussions on
the implementation of the algorithm.

\section{Proving primality in  AKS and ECPP}\label{AKSECPP}

The ECPP algorithm depends on rounds of reductions of the 
proof of primality of a  prime to the proof of primality 
of a smaller prime. The most remarkable feature of ECPP is that
a verifier who does not believe 
any conjectures can be convinced in time $\tilde{O}( \log^3 n)$
if the fast multiplication is used.
It is based on the following proposition \cite{AtkinMo93}.

\begin{proposition}
Let $N$ be an integer prime to $6$, $E$ be an elliptic curve
over $\Z/N\Z$, together with a point $P$ on $E$ and two integers
$m$ and $s$ with $s|m$. Denote the
infinite point on $E$ by $O$. For each prime divisor $q$ of $s$, 
denote $(m/q)P $ by $ (x_q:y_q:z_q)$. Assume that $mP = O$
and $gcd( z_q, N) = 1$ for all $q$. If $s > (\sqrt[4]{N} + 1)^2$,
then $N$ is a prime.
\end{proposition}

The certificate for $N$ in ECPP consists of the curve $E$,
the point $P$, $m$, $s$ and the certificate of primality of $s$.
Usually the ECPP algorithm uses elliptic curves with complex 
multiplications of small discriminants. For implementation details, 
see \cite{AtkinMo93}.

The AKS algorithm proves a number is a prime through the following proposition.

\begin{proposition}
Let $n$ be a positive integer. Let $q$ and $r$ be prime numbers.
Let $S$ be a finite set of integers. Assume 
\begin{enumerate}
\item that $q$ divides $r-1$;
\item that $n^{ r-1 \over q} \not\equiv 0, 1 \pmod{r}$;
\item that $gcd(n, b-b')=1$ for all the distinct $b, b'\in S$;
\item that ${ q + |S| - 1 \choose |S| } \geq n^{2\lfloor \sqrt{r}\rfloor}$;
\item that $(x+b)^n \equiv x^n + b  \pmod{x^r -1, n}$ for all $b \in S$. 
\end{enumerate}
Then $n$ is a power of a prime.
\end{proposition}

\section{Description and time complexity analysis of our algorithm}
\label{descriptionofalgorithm}

Now we are ready  to sketch our algorithm.

{\bf Input}: a positive integer $n$

{\bf Output}: a certificate of primality of $n$, or ``composite''. 

\begin{enumerate}
\item If $n$ is a power of an integer, return ``composite''.

\item In parallel run a composite proving algorithm,
for example, the Rabin-Miller testing \cite[Page 282]{BachSh96},  
on $n$.

\item\label{Atkinstep} 
If $n-1$ contains a prime factor between 
$\log^2 n$ and $2 \log^2 n$, skip this step.
Otherwise, call ECPP to find an elliptic curve
on $\Z/n\Z$ with $n'$ points, where
$n'$ is a  probable prime 
and $n'$ is $2$-good. Set $n = n'$.
Let $r$ be the prime factor of $n-1$ satisfying
$\log^2 n \leq r \leq 2 \log^2 n$.

\item\label{randomstep} Randomly select a number $1 < b < n$.
If $b^{n-1} \not= 1 \pmod{n}$, output ``composite''
and exit.

\item\label{not_a_pth_power} Let $a = b^{n-1 \over r^\alpha} \pmod{n}$;
If $a = 1$, or $a^{r^{\alpha-1}} = 1$,
 go back to step~\ref{randomstep}.

\item\label{extracheck}
If $gcd(a^{r^{\alpha-1}} -1, n) \not=1$,
output ``composite'' and exit.

\item\label{aksstep} 
If $( 1 + x)^n \not= 1 + x^n \pmod{n, x^r - a}, $  return ``composite'';

\item Use ECPP procedure to construct
the curve and the point and compute the order.  
Output them with $a$. Return ``prime'';
\end{enumerate}

Testing whether a number $n$ is good or not can be done in time 
$\tilde{O}(\log^3 n)$.
The step~\ref{Atkinstep} takes time $\tilde{O}(\log^{4} n)$, if 
the ECPP heuristics is true,
Conjecture~\ref{densityconjecture} in the introduction section is true,
and the fast multiplication algorithm is used.

If $n$ is indeed a prime, then the probability
of going back in  step~\ref{not_a_pth_power} is at most
$1/r$. The step~\ref{extracheck} takes
time at most $\tilde{O}(\log^2 n)$.
The step~\ref{aksstep} takes time $\tilde{O}(\log^4 n)$,
since $r \leq 2 \log^2 n$. 
Hence the heuristic expected running time of
our algorithm is $\tilde{O}(\log^{4} n)$.
Obviously the verification algorithm takes deterministic time 
$\tilde{O}(\log^4 n)$.

\section{Density of good numbers}\label{densitysection}

What is the probability that a random number has a prime factor
between $b_1$ and $b_2 = \cc b_1$? Let 
$m = \prod_{p\ \mathrm{prime,} b_1 \leq p \leq b_2} p$.
We first compute the density of integers between
$1$ and $m-1$ which has a prime factor
between $b_1$ and $b_2 $. Those numbers are precisely the
zero-divisors in $\Z/m\Z$.
The number of non-zero-divisors between $1$ and $m$ is
$ \phi(m) = m\prod_{p\ \mathrm{prime}, b_1 \leq p \leq b_2} ( 1 - {1\over p}),$
where $\phi$ is the Euler phi-function.
First we estimate the quantity:

\begin{eqnarray*}
\beta_{b_1,b_2} = 
\prod_{p\ \mathrm{prime}, b_1 \leq i \leq b_2} ( 1 - {1 \over p})\\
\end{eqnarray*}

It is known \cite{Tenenbaum95} that 
$ \prod_{p < x, p\ \mathrm{prime}} ( 1 - { 1 \over p}) 
= {e^{-\gamma} \over \ln x} ( 1 + O( {1
\over \ln x})), $
where $\gamma$ is the Euler constant.
There must exist two absolute constants $c_1, c_2$, such that
$$ {e^{-\gamma} \over \ln x} ( 1 + {c_1 
\over \ln x}) \leq \prod_{p < x, p\ \mathrm{prime}} ( 1 - { 1 \over p}) 
\leq {e^{-\gamma} \over \ln x} ( 1 + {c_2 
\over \ln x}) $$
Set $\cc = e^{c_2 - c_1 + 2}$.
\begin{eqnarray*}
\prod_{p\ \mathrm{prime}, b_1 \leq p \leq
b_2} ( 1- { 1\over p}) &=& 
{ \prod_{p\ \mathrm{prime}, p \leq
b_2} ( 1- { 1\over p}) \over \prod_{p\ \mathrm{prime}, p \leq
b_1} ( 1- { 1\over p})}\\
&\leq & { \ln b_1 \over \ln \cc b_1} 
{ 1 +  {c_2 \over \ln \cc b_1} \over 1 +  {c_1 \over \ln b_1}}\\
&=& { \ln^3 b_1 + (\ln \cc + c_2) \ln^2 b_1 \over  \ln^3 b_1 + (2 \ln \cc + c_1) 
\ln^2 b_1 + (\ln^2 \cc + 2c_1 \ln \cc) \ln b_1 + c_1 \ln^2 \cc }
\end{eqnarray*}

Thus $  1 - \beta_{b_1,b_2} = { (\ln \cc + c_1 - c_2) \ln^2 b_1 
- (\ln^2 \cc + 2c_2 \ln \cc) \ln b_1 - c_2 \ln^2 \cc
\over  \ln^3 b_1 + (2 \ln \cc + c_1) 
\ln^2 b_1 + (\ln^2 \cc + 2c_1 \ln \cc) \ln b_1 + c_1 \ln^2 \cc } > { 1  \over
\ln b_1} $, 
when $b_1$ is sufficiently large.
It is expected that the density of good primes in the set
of primes in a large interval should not be very far away from 
$ 1 \over \ln b_1$.
See Table~\ref{numericevidence} for numerical
data concerning the density of $2$-good primes around $2^{500}$.
Notice that 
\begin{eqnarray*}
\beta_{250000, 500000} &=& 0.9472455\\
1-\beta_{250000, 500000} &=& 0.0527545\\
{ 1 \over \ln 250000 } &=& 0.0804556 
\end{eqnarray*}

\begin{table}[htb!]
\caption{\label{numericevidence} \bf\large
Number of $2$-good primes around $2^{500}$}
\begin{tabular}{|c|c|c|c|c|}
\hline 
From & To & Number of primes & Number of $2$-good primes & Ratio \\
\hline
$2^{500} + 0$ &$2^{500} + 200000$ &576&35&6.07\% \\ \hline
$2^{500} + 200000$ &$2^{500} + 400000$ &558&38&6.81\% \\ \hline
$2^{500} + 400000$ &$2^{500} + 600000$ &539&30&5.56\% \\ \hline
$2^{500} + 600000$ &$2^{500} + 800000$ &568&23&4.05\% \\ \hline
$2^{500} + 800000$ &$2^{500} + 1000000$ &611&39&6.38\% \\ \hline
$2^{500} + 1000000$ &$2^{500} + 1200000$ &566&26&4.59\% \\ \hline
$2^{500} + 1200000$ &$2^{500} + 1400000$ &566&38&6.71\% \\ \hline
$2^{500} + 1400000$ &$2^{500} + 1600000$ &526&27&5.13\% \\ \hline
$2^{500} + 1600000$ &$2^{500} + 1800000$ &580&26&4.48\% \\ \hline
$2^{500} + 1800000$ &$2^{500} + 2000000$ &563&20&3.55\% \\ \hline
$2^{500} + 2000000$ &$2^{500} + 2200000$ &562&22&3.91\% \\ \hline
$2^{500} + 2200000$ &$2^{500} + 2400000$ &561&21&3.74\% \\ \hline
$2^{500} + 2400000$ &$2^{500} + 2600000$ &609&34&5.58\% \\ \hline
$2^{500} + 2600000$ &$2^{500} + 2800000$ &601&28&4.66\% \\ \hline
$2^{500} + 2800000$ &$2^{500} + 3000000$ &603&33&5.47\% \\ \hline
$2^{500} + 3000000$ &$2^{500} + 3200000$ &579&37&6.39\% \\ \hline
$2^{500} + 3200000$ &$2^{500} + 3400000$ &576&31&5.38\% \\ \hline
$2^{500} + 3400000$ &$2^{500} + 3600000$ &604&35&5.79\% \\ \hline
$2^{500} + 3600000$ &$2^{500} + 3800000$ &612&40&6.53\% \\ \hline
$2^{500} + 3800000$ &$2^{500} + 4000000$ &588&29&4.93\% \\ \hline
$2^{500} + 4000000$ &$2^{500} + 4200000$ &574&33&5.75\% \\ \hline
$2^{500} + 4200000$ &$2^{500} + 4400000$ &609&27&4.43\% \\ \hline
$2^{500} + 4400000$ &$2^{500} + 4600000$ &549&35&6.37\% \\ \hline
$2^{500} + 4600000$ &$2^{500} + 4800000$ &561&30&5.34\% \\ \hline
$2^{500} + 4800000$ &$2^{500} + 5000000$ &545&29&5.32\% \\ \hline
$2^{500} + 5000000$ &$2^{500} + 5200000$ &590&20&3.39\% \\ \hline
$2^{500} + 5200000$ &$2^{500} + 5400000$ &557&27&4.84\% \\ \hline
$2^{500} + 5400000$ &$2^{500} + 5600000$ &591&28&4.73\% \\ \hline
$2^{500} + 5600000$ &$2^{500} + 5800000$ &517&33&6.38\% \\ \hline
$2^{500} + 5800000$ &$2^{500} + 6000000$ &566&18&3.18\% \\ \hline
$2^{500} + 6000000$ &$2^{500} + 6200000$ &575&30&5.21\% \\ \hline
$2^{500} + 6200000$ &$2^{500} + 6400000$ &573&26&4.53\% \\ \hline
$2^{500} + 6400000$ &$2^{500} + 6600000$ &558&36&6.45\% \\ \hline
$2^{500} + 6600000$ &$2^{500} + 6800000$ &574&32&5.57\% \\ \hline
$2^{500} + 6800000$ &$2^{500} + 7000000$ &594&22&3.70\% \\ \hline
$2^{500} + 7000000$ &$2^{500} + 7200000$ &596&31&5.20\% \\ \hline
$2^{500} + 7200000$ &$2^{500} + 7400000$ &567&26&4.58\% \\ \hline
$2^{500} + 7400000$ &$2^{500} + 7600000$ &619&28&4.52\% \\ \hline
$2^{500} + 7600000$ &$2^{500} + 7800000$ &565&25&4.42\% \\ \hline
$2^{500} + 7800000$ &$2^{500} + 8000000$ &561&25&4.45\% \\ \hline
$2^{500} + 8000000$ &$2^{500} + 8200000$ &570&26&4.56\% \\ \hline

\end{tabular}

\end{table}

\section{Proof of the main theorem}\label{generalization}

In this section we prove the main theorem. It is built on
a series of lemmas.
Most of them are straight-forward generalizations 
of the lemmas in Berrizbeitia's paper 
\cite{Berrizbeitia02}.  We include slightly different proofs
of those lemmas, though,
for completeness. Some of the proofs are  brief, for 
details see \cite{Berrizbeitia02}.

\begin{lemma}
Let $r, p$ be primes, $r| p-1$.
If $a \in \F_p$ is not a $r$-th power of
any element in $\F_p$, then $x^r -a$ is irreducible 
over $\F_p$. 
\end{lemma}

\begin{proof}
Let $\theta$ be one of the roots of $x^r - a = 0$. 
Certainly $[\F_p (\theta) : \F_p] > 1$. 
Let $\xi \in \F_p$
be one of the $r$-th primitive roots of unity.
$$ x^r - a = x^r - \theta^r = \prod_{0\leq i\leq r-1} (x-\xi^i \theta). $$
Let $[\F_p (\theta) : \F_p] = r'$. 
Then for all $i$, $[\F_p(\xi^i
\theta):\F_p] = r'$. Hence $x^r - a $ will be factored into
polynomials of degree $r'$ only. 
Since $r$ is a prime,
this is impossible, unless that $r' = r$.
\end{proof}

\begin{lemma}
Let $n > 2$ be an integer. Let $r$ be a prime and $r^\alpha || n-1$.
Suppose that there exists a integer $ 1 < a < n$ such that
\begin{enumerate}
\item $a^{r^\alpha} \equiv 1 \pmod{n}$;
\item $gcd(a^{r^{\alpha-1}} - 1, n) = 1$;
\end{enumerate}
Then there must exist a prime factor $p$ of $n$, such
that $r^\alpha || p-1$ and
$a$ is not a $r$-th power of any element in $\F_p$.
\end{lemma}

\begin{proof}
For any prime factor $q$ of $n$, $ a^{r^\alpha} \equiv 1 \pmod{q}$ 
and $ a^{r^{\alpha-1}}  \not\equiv 1 \pmod{q}$, so $ r^\alpha | q - 1$. 
If $ r^{\alpha+1} | q -1 $ for all
the prime factors, then $r^{\alpha+1} | n -1$, contradiction.
Hence there exists a prime factor $p$, such that $r^\alpha || p-1$.
Let $g$ be a generator in $\F_p^*$. If $a = g^t $ in $\F_p$,
then $p-1 | tr^\alpha$, and $p-1 \not| tr^{\alpha-1}$. Hence $r \not| t$.
\end{proof}

In the following text, we assume that $n$ is an integer,
$n = p^l d$ where $p$ is a prime and $gcd(p,d)=1$.
Assume $r$ is a prime and $r|p-1$.
Let $x^r -a$ be an irreducible polynomial in $\F_p$. 
Let $\theta$ be one of the roots of $x^r - a$. 
For any element in the field $\F_p (\theta)$, we can find a
unique polynomial $f\in \F_p [x]$ of degree less than $r$
such that the element can be represented by $f(\theta)$.
Define
$\sigma_m: \F_p (\theta) \rightarrow \F_p (\theta)$
as $\sigma ( f(\theta) ) = f(\theta^m)$.

\begin{lemma}
We have that $a^m = a$ in $\F_p$ iff $ \sigma_m \in Gal(\F_p(\theta)/\F_p)$.
\end{lemma}

\begin{proof}
($\Leftarrow$): Since $ \sigma_m \in Gal(\F_p(\theta)/\F_p)$,
$\theta^m$ must be a root of $x^r - a$. Hence $a= (\theta^m)^r = a^m $ in
$\F_p$.

($\Rightarrow$): For any two elements $a,b\in \F_p (\theta)$, we need to
prove that $\sigma_m (a+b)  = \sigma_m (a) + \sigma_m (b) $ and
$\sigma_m (ab)  = \sigma_m (a) \sigma_m (b) $. The first one is
trivial from the definition of $\sigma_m$. Let $a =f_a (\theta)$ and
$b = f_b(\theta)$ where $f_a (x), f_b(x) \in \F_p [x]$ has degree less
than $r-1$. If $ deg(f_a(x) f_b(x)) \leq r - 1$, 
it is easy to see that $\sigma_m (ab)  = \sigma_m (a) \sigma_m (b) $.
Now assume that $ deg(f_a(x) f_b(x)) \geq r $.
Then $f_a (x) f_b (x) = h(x) + (x^r - a) p(x)$ where 
$h(x), p(x) \in \F_p[x]$ and $deg (h(x)) < r$.
Then $\sigma_m (ab) = \sigma_m ( h(\theta)) = h(\theta^m)
= h(\theta^m) + (a^m - a) p(\theta^m)
= h(\theta^m ) +  (\theta^{mr} - a) p(\theta^m) = f_a (\theta^m)
f_b(\theta^m) = \sigma_m (a) \sigma_m (b)$.

This shows that $\sigma_m$ is a homomorphism. Now we need to
prove that it is one-to-one. It is obvious since
$\theta^m$ is a root of $x^r - a=0$.
\end{proof}

Define $G_m = \{f(\theta) \in \F_p(\theta)^* | f(\theta^m ) = f(\theta)^m\}$.
It can be verified that $G_m$ is a group when $\sigma_m$ 
is in $Gal(\F_p(\theta)/\F_p)$.

\begin{lemma}\label{lemmaGn}
Suppose $\sigma_n \in Gal(\F_p(\theta)/\F_p)$. Then for any 
$i,j\geq 0$, $\sigma_{d^i p^j} \in Gal(\F_p(\theta)/\F_p)$
and $G_n \subseteq G_{d^i p^j}$.
\end{lemma}

\begin{proof}
Notice that the map $x\rightarrow x^{p^l}$ is 
a one-to-one map in $\F_p (\theta)$. The equation $a^n =  a$ 
implies that $(a^d)^{p^l} = a$, hence $a^d = a$, and $a^{d^i p^j} =a$.
We have $\sigma_{d^i p^j} \in Gal(\F_p(\theta)/\F_p)$.

Let $f(\theta) \in G_n$. Thus $ f( \theta^n ) = f(\theta)^n$, this implies
$f(\theta^{p^ld}) = f(\theta)^{p^ld} = f(\theta^{p^l})^d$. So $\theta^{p^l}$ is
a solution of $f(x^d) = f(x)^d$. Since it is 
one of the conjugates of $\theta$, $\theta$ must be a solution as well.
This proves that $f(\theta^d ) = f(\theta)^d$. Similarly
since $\theta^d$ is also one of the conjugates of $\theta$,
 as $\sigma_d \in Gal(\F_p(\theta)/\F_p)$,  we have
$f(\theta^{d^2} ) = f(\theta^d)^d
=f(\theta)^{d^2}$. By reduction,
$f(\theta^{d^i} ) = f(\theta)^{d^i}$ for $k\geq 0$.
Hence $f(\theta^{d^ip^j} ) = f(\theta^{d^i})^{p^j}
= f(\theta)^{d^i p^j}$. This implies that 
$f(\theta) \in G_{d^ip^j}$.
\end{proof}

\begin{lemma}
If $\sigma_{m_1}, \sigma_{m_2} \in Gal(\F_p(\theta)/\F_p)$
and $\sigma_{m_1} = \sigma_{m_2}$, then
$|G_{m_1} \cap G_{m_2}| $ divides $m_1 - m_2$.
\end{lemma}

This lemma is straight forward from the definition.

\begin{lemma}
Let $A = a^{r^{\alpha-1}}$.
If $(1 + \theta) \in G_n$, so is $ 1+ A^i  \theta$ for 
any $i =1,2,3,\cdots, r-1$. And $|G_n| \geq 2^r$.
\end{lemma}

\begin{proof}
If $(1 + \theta) \in G_n$, this means that
$ ( 1+ \theta)^n = 1 + \theta^n$. It implies that
 $( 1+ \theta')^n = 1+\theta'^n$
for  any conjugate $\theta'$  of $\theta$.
Since  $A$ is a primitive root of unity in $\F_p$,
hence $A^i \theta$ are conjugates of $\theta$. We have
$ ( 1 + A^i \theta)^n = 1 + (A^i \theta)^n = 1 + (A^n)^i \theta^n$
and we know that $A^n = A$. This proves that $1+ A^i \theta \in G_n$.
The group $G_n$ contains all the elements in the set
$$ \{ \prod_{i=0}^{r-1} ( 1 + A^i \theta)^{\epsilon_i} |
\sum_{i=0}^{r-1} \epsilon_i < r \}, $$
by simple counting we have $|G_n| \geq 2^r$.

\end{proof}

Finally we are ready to give the proof of 
the main theorem  (Theorem~\ref{main}) of this paper.

\begin{proof}
Since $|Gal(\F_p(\theta)/\F_p)| = r$, hence there exist two different
pairs $(i_1,j_1) $ and $(i_2,j_2) $
with $0\leq i_1,j_1, i_2, j_2\leq \lfloor \sqrt{r} \rfloor$, 
such that $\sigma_{d^{i_1}p^{j_1}} = \sigma_{d^{i_2}p^{j_2}}$.
According to Lemma~\ref{lemmaGn}, $G_n \subseteq G_{d^{i_1}p^{j_1}}$,
$G_n \subseteq G_{d^{i_2}p^{j_2}}$, this implies that
$G_n \subseteq G_{d^{i_1}p^{j_1}} \cap G_{d^{i_2}p^{j_2}}$.
Therefore $|G_n|$ divides $d^{i_1}p^{j_1} - d^{i_2}p^{j_2}$,
but $d^{i_1}p^{j_1} - d^{i_2}p^{j_2} < n^{\lfloor \sqrt{r} \rfloor}
\leq 2^{ \sqrt{r} \log n} \leq 2^r$.
hence $d^{i_1}p^{j_1} - d^{i_2}p^{j_2} = 0$,
which in turn implies that $n $ is a power of $p$. 
\end{proof}

\section{Implementation and conclusion}\label{implementation}

In this paper, we propose a random primality proving algorithm
which runs in heuristic time $\tilde{O}(\log^4 n)$. It generates 
a certificate of primality of length $O(\log n)$ 
which can be verified in deterministic time
$\tilde{O}(\log^4 n)$.

When it comes to implement the algorithm, space is a   bigger issues 
than time. Assume that $n$ has $1000$ bit, which is
the range of practical interests.
To compute $(1 +x)^n \pmod{n, x^r-a}$,
we will have an intermediate polynomial of size $2^{30}$ bit, or 
$128$M bytes. As a comparison, ECPP is not very demanding on space.
In order to make the algorithm available on a desktop PC,
space efficient exponentiation of $1+x$ is highly desirable.
This is the case  for the original version of the 
AKS algorithm as well.

For the sake of theoretical clarity, we use just one round
of ECPP reduction in the algorithm. To implement the algorithm,
it may be better to follow the
ECPP algorithm and launch the iteration of AKS
as soon as an intermediate prime becomes good.
Again assuming that the intermediate primes are 
distributed randomly in the range, the expected 
number of rounds will be $\log \log n$. 
It is a better strategy since the intermediate
primes get smaller and smaller.

\paragraph{Acknowledgements:} We thank Professor Pedro Berrizbeitia
for very helpful discussions and comments.

\bibliographystyle{plain}
\bibliography{crypto}

\end{document}